\documentclass{amsart}
\usepackage{epsfig}
\usepackage{amscd}

\newtheorem{theorem}{Theorem}
\newtheorem{lemma}{Lemma}
\newtheorem{cor}{Corollary}
\newtheorem{prop}{Proposition}

\newtheorem{ex}{Example}

\newcommand{\tr}{{\rm tr}}

\newcommand{\cL}{{\mathcal L}}
\newcommand{\cR}{{\mathcal R}}

\newcommand{\cI}{{\mathcal I}}
\newcommand{\cJ}{{\mathcal J}}

\newcommand{\cW}{{\mathcal W}}
\newcommand{\cS}{{\mathcal S}}
\newcommand{\cT}{{\mathcal T}}
\newcommand{\bC}{\hspace{-0.3pt}{\mathbb C}\hspace{0.3pt}}
\newcommand{\g}{SL_2(\bC)} 
\newcommand{\ch}{\chi_\rho}

\begin{document}

\title{Rings of $\g$-Characters and the Kauffman Bracket Skein Module}
\author{Doug Bullock \\ Boise State University, Boise, ID 83725, USA}

\maketitle

\begin{abstract} Let $M$ be a compact orientable 3-manifold.
The set of characters of $\g$ representations of the fundamental group
of $M$ forms a closed affine algebraic set.  We show that its
coordinate ring is isomorphic to a specialization of the Kauffman
bracket skein module modulo its nilradical. This is accomplished
by making the module into a combinatorial analog of the ring, in which
tools of skein theory are exploited to illuminate relations among
characters.  We conclude with an application, proving that a small
manifold's specialized module is necessarily finite dimensional.
\par\mbox{}\par\noindent {\em Keywords: knot, link, skein theory,
representation theory, 3-manifold.}  \par\mbox{}\par\noindent {\em AMS
(MOS) Subject Classification:} 57M99.
\end{abstract}

\section{Introduction}

The Kauffman bracket skein module is an invariant of 3-manifolds
which, until recently, was both difficult to compute and topologically
mysterious. The discovery \cite{estimate} that a specialization of the
module dominates the ring of $\g$ characters of the fundamental group
shed some light on the meaning of the module.  The relationship also
provided estimates \cite{reps} \cite{estimate} of the module's size
and computational tools \cite{BP}.  The central result of this paper
sharpens the focus considerably, for we show that the specialization,
modulo its nilradical, is exactly the ring of characters.

The construction depends upon a natural correspondence between knots
and functions on the set of characters.  Given an orientation, a knot
corresponds to a conjugacy class in the fundamental group of a
3-manifold $M$.  A formal linear combination of knots is therefore a
template upon which one may evaluate characters of $\pi_1(M)$
represented in $\g$.  Conversely, one may interpret polynomial
functions on the set of characters as linear combinations of links.
These functions form an algebra, of which the the ring of characters
is a quotient.  The skein module is also a quotient (of the linear
space of links) and it, too, has a ring structure.  The correspondence
between knots and functions descends to these quotients, where its
kernel is exactly the nilpotent elements of the skein module.

The proof proceeds in three stages, the first of which (Section 2)
recapitulates parts of \cite{estimate} and \cite{CS}.  We cover the
necessary background, including precise definitions of the principal
objects of study.  Once the vocabulary is in place we define the map,
$\Phi$, taking knots to functions on the character set.  The proof
that it descends to a surjection on the quotients is quite simple,
depending primarily on the following observation: the Kauffman bracket
skein relation maps to the fundamental $\g$ trace identity,
\[\tr(AB)+\tr(AB^{-1})=\tr(A)\tr(B).\]
The characterization of $\ker \Phi$, however, is significantly more involved.

We first study the correspondence for handlebodies and free groups,
beginning in Section 3 with an investigation of trace identities.  We
define a map, $\Psi$, which sends functions on characters into the
skein module.  It turns out that $\ker \Phi$ contains exactly those
trace identities which $\Psi$ does not send to zero.  The main result
recalls work of Procesi \cite{procesi} and Razmyslov \cite{raz}, who
classified a large family of trace identities on arbitrary matrix
rings.  We show that these identities, when restricted to $\g$, map to
zero.  The skein module emerges here as a useful combinatorial tool.
Although it is now possible to give purely algebraic proofs of the
results in this section, they were all discovered by experimenting
with the skein module.  We have retained the geometric arguments for
they serve to illustrate the interplay between the two theories
and---once the reader has become comfortable with the skein
moves---they make for shorter proofs.

Section 3 attains sufficient conditions for a trace identity to vanish
in the skein module; Section 4 provides the finishing touch.  We rely
on a defining set of polynomials for the character set given in
\cite{GM}, the central results of which are reiterated in an effort to
keep this paper self contained.  Most of these polynomials turn out to
be specialized Procesi identities, while the remaining few succumb to
other tools from Section 3.  It follows from a standard result of
algebraic geometry that the only trace identities not vanishing in the
skein module are nilpotent.  It is then a small step to extend the
result to arbitrary compact 3-manifolds.

The author would like to thank Professors Charles Frohman, Xaio-Song
Lin, J\'{o}zef Przytycki and Bruce Westbury for many helpful
conversations and suggestions; Adam Sikora in particular for his
insight into the importance of nilpotents; and the organizers and
participants of the Banach Center's Mini Semester on Knot Theory,
where the ideas in this paper first began to coalesce.

\section{Definitions and Background}

Let $M$ be a compact orientable 3-manifold.  
The Kauffman bracket skein module of $M$ is
an algebraic invariant, denoted $K(M)$, which is built from the set $\cL_M$ of framed
links in $M$.  By a framed link we mean an embedded collection of
annuli considered up to isotopy in $M$, and we include the empty
collection $\emptyset$.  Three links $L$, $L_0$ and $L_\infty$ are
said to be {\em Kauffman skein related} if they can be embedded
identically except in a ball where they appear as shown in Figure 1
(framings are vertical with respect to the page).
  \begin{figure}
    \mbox{}\hfill\epsfig{figure=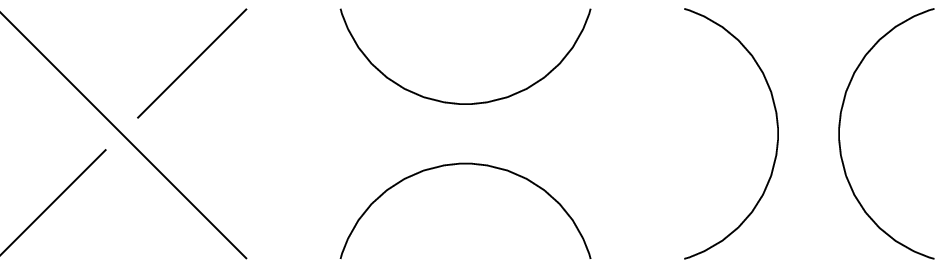,width=2in}\hfill\mbox{}
    \newline\mbox{}\hspace{1.7in}\makebox[1.65in]{$L$\hfill$L_0$\hfill
       $L_\infty$}
    \caption{}
  \end{figure}
The notation $L \amalg \bigcirc$ indicates the union of $L$ with an
unlinked 0-framed unknot. 
      
Let $R$ denote the ring of Laurent polynomials $\bC[A^{\pm 1}]$ and
$R{\mathcal L}_{M}$   the free $R$-module with basis ${\mathcal
L}_{M}$. If $L$, $L_0$ and $L_\infty$ are Kauffman skein related then
$L-AL_{0}-A^{-1}L_{\infty}$ is called a {\em skein relation}.  For any
$L$ in ${\mathcal L}_M$ the expression $L \amalg \bigcirc + (A^2
+A^{-2})L$ is called a {\em framing relation}.  Let $S(M)$ be the
smallest submodule of $R{\mathcal L}_M$ containing all possible skein
and framing relations.  We define $K(M)$ to be the quotient
$R{\mathcal L}_{M}/S(M)$.

The indeterminate $A$ is often interpreted as a complex number so that
$K(M)$ becomes a vector space.  It seems that the simplest value is
$A=-1$, and we let $V(M)$ denote this specialization.  Notice that the
specialized skein relations imply \par\mbox{}\par\noindent
\mbox{}\hfill\epsfig{figure=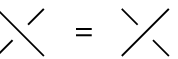}\hfill\mbox{}
\par\mbox{}\par\noindent in $V(M)$.  There is a product on links,
$L_1L_2=L_1\cup L_2$, which makes $V(M)$ into a commutative algebra
with $\emptyset$ serving as the identity.  If follows from
\cite[Theorem 1]{states} that $V(M)$ is generated by a finite set of
knots.

By a {\em representation}  we mean a homomorphism of
groups 
\[\rho : \pi_1(M) \rightarrow \g.\]  
The {\em character} of a representation is the composition
\[\ch = {\rm trace}\circ\rho,\]
and $X(M)$ denotes the set of all characters.
For each $\gamma \in \pi_1(M)$ there is a function $t_\gamma :
X(M) \rightarrow \bC$ given by $\ch \mapsto \ch(\gamma)$.  The
following theorem appears to have been discovered independently 
by Vogt \cite{V} and Fricke \cite{F}, first proved by Horowitz
\cite{H}, and then rediscovered by Culler and Shalen \cite{CS}. 

\begin{theorem} {\rm (Vogt, Fricke, Horowitz, Culler--Shalen)}
There exits a finite set of elements $\{\gamma_1,\ldots,\gamma_m\}$ in $\pi_1(M)$ such that every
$t_\gamma$ is an element of the polynomial ring $\bC[t_{\gamma_1},\ldots,t_{\gamma_m}]$.  
\end{theorem}

For Culler and Shalen, Theorem 1 was an initial step in a much deeper result.

\begin{theorem} {\rm (Culler--Shalen)} If every $t_\gamma$ is an
element of $\bC[t_{\gamma_1},\ldots,t_{\gamma_m}]$, then $X(M)$ is a
closed algebraic subset of $\bC^m$.  
\end{theorem}

Recall that a {\em closed algebraic set} $X$ in $\bC^m$ is the
common zero set of an ideal of polynomials in $\bC[x_1,\ldots,x_m]$.
The elements of  $\bC[x_1,\ldots,x_m]$ are {\em polynomial functions} on
$X$, and the functions $x_i$ are {\em coordinates} on $X$.  The quotient
of $\bC[x_1,\ldots,x_m]$ by the ideal of polynomials vanishing
on $X$ is called the {\em coordinate ring} of $X$.  Different choices
of coordinates would clearly lead to different parameterizations of
$X$, but it follows from \cite{CS} that any two parameterizations of
$X(M)$ are equivalent via polynomial maps.  Hence
their coordinate rings are isomorphic and we may identify them as one
object: the {\em ring of characters} of  $\pi_1(M)$, which we denote
by $\cR(M)$.  

Each knot $K$ determines a unique $t_\gamma$ as follows.
Let $\vec{K}$ denote an unspecified orientation on $K$. Choose any
$\gamma \in \pi_1(M)$ such that $\gamma \simeq \vec{K}$ (meaning the
loop $\gamma$ is freely homotopic to an embedding of $\vec{K}$).   
Since trace is invariant under conjugation it makes sense to define
$\ch(\vec{K})= \ch(\gamma)$.  Since $\tr(A)
= \tr(A^{-1})$ in $\g$ we can also define $\ch(K)=\ch(\gamma)$.  Thus
$K$ determines the map $t_\gamma$. Conversely, any $t_\gamma$ is
determined by some (non-unique) $K$. The main theorem of
\cite{estimate} is that this correspondence is well  defined at the
level of $V(M)$. 

\begin{theorem} The map $\Phi : V(M) \rightarrow \cR(M)$ given by 
\[ \Phi(K)(\ch) = -\ch(K) \]
is a well defined surjective map of algebras.  If $V(M)$ is generated
by the  knots $K_1,\ldots,K_m$ then $-\Phi(K_1),\ldots,-\Phi(K_m)$
are coordinates on $X(M)$.
\end{theorem}

\begin{proof}
Let $\bC^{X(M)}$ denote the algebra of
functions from  $X(M)$ to $\bC$.  Define a map 
\[ \widetilde{\Phi} : \bC\cL_M \rightarrow \bC^{X(M)} \]
as follows.  If $K$ is a knot set
\[ \widetilde{\Phi}(K)(\ch) = -\ch(K). \]
If $L$ is a link with components $K_1,\ldots,K_n$ set
\[ \widetilde{\Phi}(L) = \prod_{i=1}^n \widetilde{\Phi}(K_i). \]
Set $ \widetilde{\Phi}(\emptyset) = 1$ and extend linearly. 

Consider the image of $S(M)$ under $\widetilde{\Phi}$.  For a framing
relation, $L \amalg \bigcirc + 2L$, we have
\begin{equation*}\begin{split}
\widetilde{\Phi}(L \amalg \bigcirc + 2L)(\ch) & = \widetilde{\Phi}(L)
            \widetilde{\Phi}(\bigcirc + 2\;\emptyset)\\
& =  -\ch(\bigcirc) + 2\\
& =  -\tr({\rm Id}) + 2\\
& =  0.
\end{split}\end{equation*}
Next, let  $L + L_0 + L_\infty$ be a skein relation in which $L$ and $L_0$
are knots.  It follows that $L_\infty$ has two components, $K_1$ and
$K_2$.  Assume embeddings as in Figure 1 and choose
a base point $*$ in the neighborhood where $L$, $L_0$ and $L_\infty$
differ.  It is now possible to find loops $a$ and $b$ in $\pi_1(M,*)$
so that a slight perturbation of $ab$ gives $\vec{L}$.  With favorable
orientations on the other knots we
have $ab^{-1}\simeq\vec{L}_0$, $a \simeq \vec{K}_1$, and $b
\simeq\vec{K}_2$.  Given any $\ch$,  set
$A=\rho(a)$ and $B = \rho(b)$ so that
\begin{equation*}\begin{split}
\widetilde{\Phi}(L+L_0+L_\infty)(\ch)&=-\ch(L)-\ch(L_0)+\ch(K_1)\ch(K_2)\\
& = -\tr(AB)-\tr(AB^{-1})+\tr(A)\tr(B) \\
& = 0.
\end{split}\end{equation*}
Finally, note that every skein relation can be written as $L' \cup L +
L' \cup L_0 + L' \cup L_\infty$ where $L$ and $L_0$ are knots.  Hence
$\widetilde{\Phi}$ descends to a well defined map of algebras,
\[ \Phi : V(M) \rightarrow \bC^{X(M)},\]
which is determined by its values on knots.  

Let $K_1, \ldots, K_m$ be generators of $V(M)$. Every element of $V(M)$
can be written as a polynomial in these knots, so the image of $\Phi$
lies in $\bC[-\Phi(K_1),\ldots,-\Phi(K_m)]$.  Since
each $t_\gamma$ is equal to $-\Phi(K)$ for some knot $K$, Theorems 1
and 2 imply that
the functions $-\Phi(K_i)$ are coordinates on $X(M)$.  It follows
that $\Phi$ maps onto $\cR(M)$.  
\end{proof}

\section{Trace Identities}

In the previous section we obtained a surjection $\Phi : V(M)
\rightarrow \cR(M)$ based on a natural correspondence between knots
and functions on $X(M)$.  Under this correspondence elements of $S(M)$
were sent to polynomials that vanish on $X(M)$, making $\Phi$ well
defined.  Our ultimate goal is to show that $\ker \Phi$ is the set of
nilpotent elements in $V(M)$.  To this end we reverse the
correspondence, mapping polynomials on $X(M)$ to elements of $V(M)$.
For now, we will treat only the case where $M$ is a handlebody.  In
this setting the kernel of $\Phi$ consists of polynomials that vanish
on $X(M)$ but not in $V(M)$.

For the time being we will be concerned only with free groups, so
throughout this section $H$
will denote the manifold $P \times I$ where $P$ is the planar
surface in Figure 2.  
 \begin{figure}
    \mbox{}\hfill\epsfig{figure=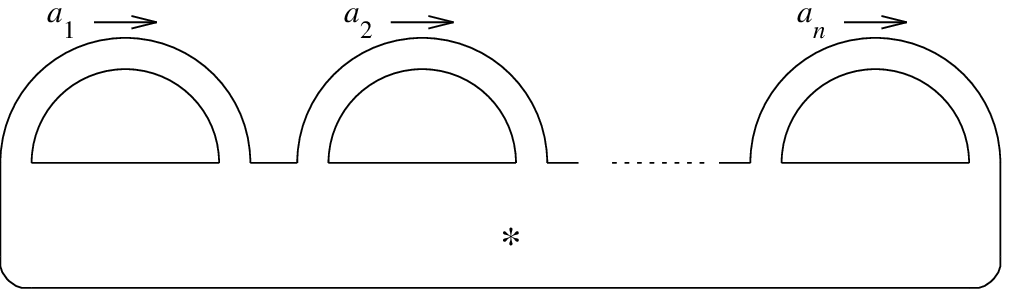,width=4in}\hfill\mbox{}
    \caption{}
 \end{figure}
We also fix a base point $*$ in $P$ and a set of
generators $\{a_1,\ldots,a_n\}$ for $\pi_1(H,*)$.  Each loop $a_i$
travels once across the $i$-th handle in the direction shown in Figure
2. Let $\cW$ denote $\pi_1(H,*)$ modulo the
equivalence
\[ w \sim w' \Longleftrightarrow \mbox{$w'=w^{-1}$ or $w' = gwg^{-1}$
for some $g\in\pi_1(H,*)$.}\]
Consider the ring of polynomials $\bC[\cW]$.  

\begin{ex}
\[ p = (a_1)(a_2)(a_3) - (a_1a_2)(a_3)-(a_1a_3)(a_2)-(a_2a_3)(a_1) +
(a_1a_2a_3) + ( a_1a_3a_2) \]
\end{ex}

\begin{ex}
\begin{equation*}\begin{split}
q& =    (a_1)^2+(a_2)^2+(a_3)^2+(a_1a_2)^2+(a_1a_3)^2+(a_2a_3)^2\\
 &\quad+(a_1a_2a_3)^2+(a_1a_2)(a_1a_3)(a_2a_3)
        +(a_1a_2a_3)(a_1)(a_2)(a_3)\\
 &\quad-(a_1a_2a_3)(a_1)(a_2a_3)-(a_1a_2a_3)(a_2)(a_1a_3)
        -(a_1a_2a_3)(a_3)(a_1a_2)\\
 &\quad-(a_1)(a_2)(a_1a_2)-(a_1)(a_3)(a_1a_3)-(a_2)(a_3)(a_2a_3)-4
\end{split}\end{equation*}
\end{ex}

The parentheses are necessary to distinguish multiplication in
$\pi_1(H)$ from multiplication in $\bC[\cW]$. Note that there is some
ambiguity in the notation for an individual element of $\bC[\cW]$.
For instance $(w^2)+(1)-(w)^2$ is the same as
$(ww)+(ww^{-1})-(w)(w^{-1})$.  Occasionally it will be convenient
to write a polynomial using non-reduced words.

A representation of $\pi_1(H,*)$ in $\g$ is any assignment of matrices to each
$a_i$.  Letting parentheses denote the operation of trace, each
element of $\bC[\cW]$ becomes a function from the representation space to $\bC$.  The elements of $\bC[\cW]$ that
vanish as functions on the set of representations are called $\g$ {\em
trace identities}.  They form an ideal $\cI \subset \bC[\cW]$. 
 
Each $w \in \cW$ corresponds to a unique unoriented curve up to
homotopy in $H$.  We will use $K_w$ to denote any knot in this
homotopy class.  Since crossings are irrelevant, $K_w$ represents a
unique element of $V(H)$.  The assignment $w \mapsto -K_w$ defines
a surjection of algebras,
\[ \Psi : \bC[\cW] \rightarrow V(H).\]

The map $\Psi$ turns an element of $\bC[\cW]$ into
a linear combination of links in $V(H)$, where we can apply skein
theory.  The basic tool for calculating in  $V(H)$ is a resolving tree.   Let
$T$ be a finite, connected, contractible graph in which no vertex has
valence greater than three.  Assume that each vertex is labeled $cL$ 
for some $c \in \bC$ and some $L \in \cL_H$.  Assume further that there is a
distinguished vertex $c_0L_0$
called the {\em root}.  Define the potential of a vertex to be the
number of edges in a path to the root.   A (necessarily
univalent) vertex that is not adjacent to one of higher potential is called a
{\em leaf}.    We say $T$ is a {\em resolving tree} for $c_0L_0$ if each
vertex $cL$ satisfies exactly one of the following.
\begin{enumerate}
\item  $cL$ is a leaf.
\item  $cL$ is adjacent to exactly one higher potential vertex, $-2cL$.
\item  $cL$ is adjacent to exactly two higher potential vertices, $c'L'$ and
$c''L''$, in which case  $cL-c'L'-c''L''$ is a framing relation.
\end{enumerate}

\begin{figure}[b]
\epsfig{figure=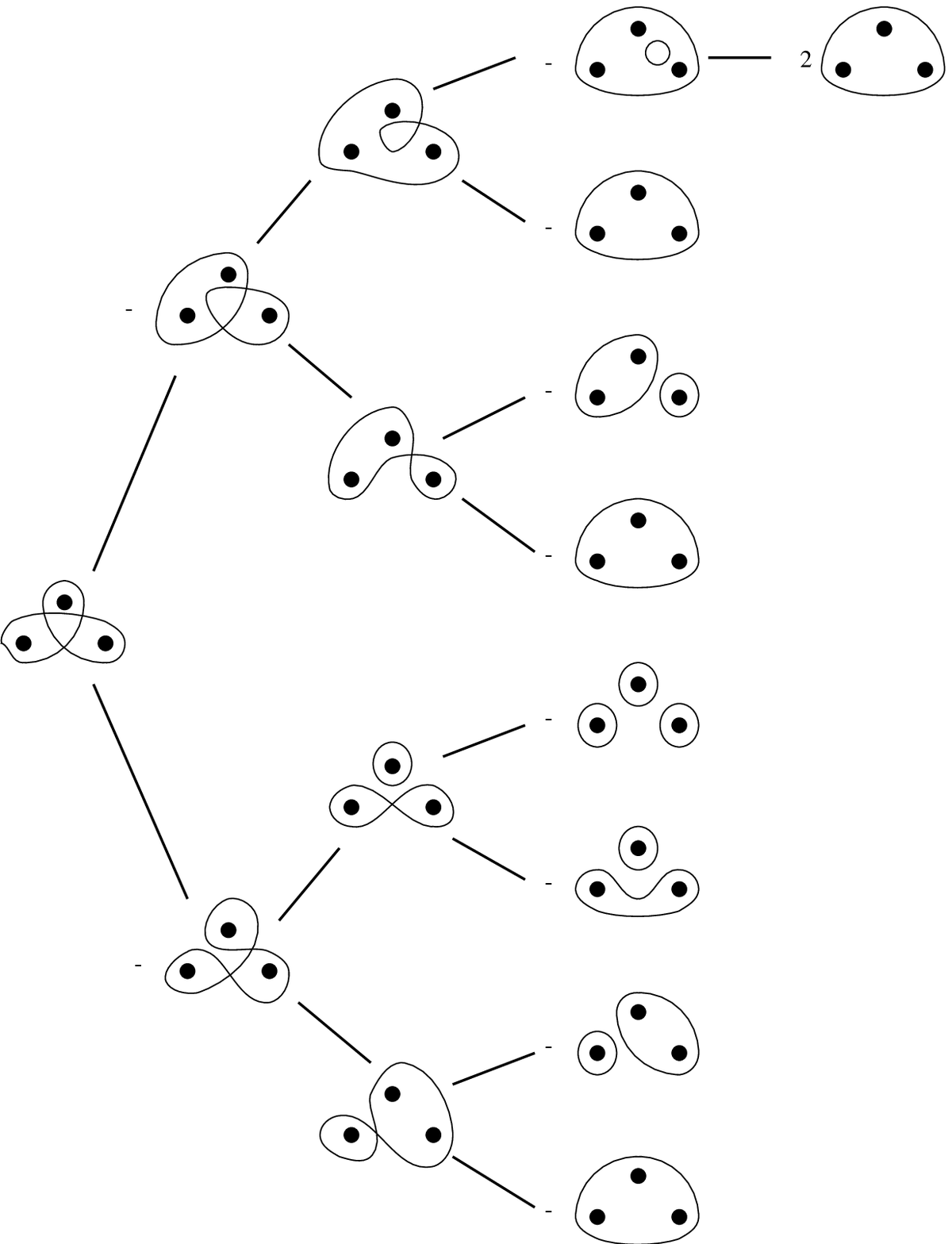,width=5in} \caption{}
\end{figure}

Figure 3, in which the dots represent a thrice punctured plane, is an
example of a resolving tree for any knot that projects to the
leftmost diagram.  It is also an example of the most common way to 
produce a resolving tree.  Beginning with a projection of the root,
the tree grows by smoothing one crossing at a time.  Once all
crossings have been eliminated, trivial circles are removed via
framing relations.  Summing over the leaves gives the {\em standard
resolution} of the root---an element of $\bC\cL_H$ which is equal
to the root in $V(H)$.  Although the procedure given here does not
result in a unique tree, the following theorem \cite{skeinmodules}
implies uniqueness of the standard resolution in $\bC\cL_H$.

\begin{theorem} {\rm (Przytycki)} The links in $H$ represented by
diagrams in $P$ with no crossings and no trivial circles are a basis
for $V(H)$.
\end{theorem} 

A {\em resolving forest} for an element of
$\bC\cL_H$ is simply a collection of trees, one for each term in the
linear combination.  As with individual links, there is a standard
resolution of each element of $\bC\cL_H$.   Summing the potential
function over all vertices assigns a useful complexity to a forest,
the {\em total potential}.  

The remainder of this section is devoted to the establishment
of conditions under which $\Psi$ maps an identity to zero.

\begin{lemma}
In Examples 1 and 2 we have $\Psi(p)=\Psi(q)=0$.
\end{lemma}

\begin{proof}
The root and leaves of the tree in Figure 3 are links
representing the terms of $p$.  One may check that the identity in
$V(H)$ given by this resolution is precisely $\Psi(p)$.  
For $q$, resolve the diagram in Figure 4,
  \begin{figure}
    \mbox{}\hfill\epsfig{figure=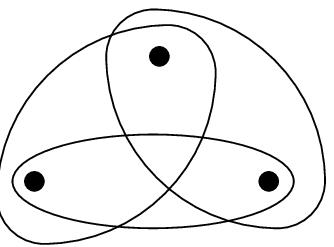,width=1in}\hfill\mbox{}
    \caption{}
  \end{figure} 
which represents $-(a_1a_2)(a_1a_3)(a_2a_3)$.
\end{proof}

If $p$ is a trace identity then a natural way to produce a new trace
identity, $q$, is to substitute new words for each $a_i$ in $p$.  If
$\Psi(p)=0$ then one would hope $\Psi(q)=0$ as well.  Although this is
true, the proof requires some effort.  

\begin{lemma}
Let $p \in \bC[\cW]$.  If there exist words $w_1$ and $w_2$ such that
$(w_1w_2)+(w_1w_2^{-1})-(w_1)(w_2)$ divides $p$ then $\Psi(p)=0$.
Also, if $p$ is divisible by $(1)-2$ then $\Psi(p)=0$. 
\end{lemma}

\begin{proof}
In the first case consider the loop $w_1w_2$, but perturbed slightly
so as to become an embedding.  By definition the resulting knot is
some $K_{w_1w_2}$.  Similarly, perturb $w_1w_2^{-1}$, $w_1$ and $w_2$
to obtain embeddings of $K_{w_1w_2^{-1}}$, $K_{w_1}$ and $K_{w_2}$.
The perturbations may be chosen so that the embeddings of
$K_{w_1w_2}$, $K_{w_1w_2^{-1}}$ and $K_{w_1}K_{w_2}$ coincide outside
of a small neighborhood of the base point.  Within that neighborhood
they appear as in Figure 1, so they form a Kauffman skein triple.  We
now have
\begin{equation*}
\begin{split}
0 & = -K_{w_1w_2}-K_{w_1w_2^{-1}}-K_{w_1}K_{w_2}\\
  & = \Psi((w_1w_2)+(w_1w_2^{-1})- (w_1)(w_2)),
\end{split}
\end{equation*}
which implies $\Psi(p)=0$.  In the second case $\Psi(p)$ contains a
factor of $\bigcirc + 2\;\emptyset$, which also implies $\Psi(p)=0$.  
\end{proof}

\begin{prop}
Let $p \in \bC[\cW]$.  Choose words $w_1,\ldots,w_n$, and form a
new polynomial $q$ by substituting $w_i$ for $a_i$ in $p$.  If
$\Psi(p)=0$ then $\Psi(q)=0$.
\end{prop} 

\begin{proof} The proof is by induction on a complexity, $\kappa(p)$,
which we define as follows.  For each $w \in \cW$ choose a diagram in
$P$ representing $K_w$.  Express $\Psi(p)$ as an element of $\bC\cL_H$
using these diagrams, and then 
choose a forest for its standard resolution. Define $\kappa(p)$ to be
the minimum total potential over all choices of diagrams and
forests. 

Assume first that $\kappa(p)=0$, implying diagrams in which 
$\Psi(p)$ is expressed as its own standard
resolution.  If $\Psi(p)=0$, we can invoke Theorem 4 to
conclude that this particular expression of $\Psi(p)$ is formally zero
in $\bC\cL_H$.  It is not possible for a diagram to represent more
than one $w$, so $p$ (and hence $q$) must be identically zero.  

Now assume that $\kappa(p)>0$.   Choose diagrams and a forest realizing
$\kappa(p)$; also select a root $cL$ which  is not a leaf.  There are
three cases depending on the first resolution of $cL$. 

{\em Case 1: The resolution removes a self crossing of some
component.}  Letting $K$ denote that component, we construct loops in
$\pi_1(H,*)$.  Begin by choosing a point $x$ near the
crossing in question.  Let $\alpha_0$ be an arc running from
$*$ to $x$; let $\alpha_1$ be an arc running parallel to $\vec K$
until it returns to $x$; and let $\alpha_2$ be an arc parallel
to the remaining  portion of $\vec K$. Set $\gamma_1 =
\alpha_0\alpha_1\alpha_0^{-1}$ and $\gamma_2 =
\alpha_0\alpha_2\alpha_0^{-1}$.  
We now have $K = K_{\gamma_1\gamma_2}$.
Furthermore, the resolution changes $K$ into $K_{\gamma_1\gamma_2^{-1}}$ and
$K_{\gamma_1}K_{\gamma_2}$.

The term of $p$ represented by $cL$ must contain the indeterminate
$(\gamma_1\gamma_2)$.  Replace that appearance
of $(\gamma_1\gamma_2)$ with $(\gamma_1)(\gamma_2)-(\gamma_1\gamma_2^{-1})$, creating a new
polynomial $p'$.  Since $p-p'$ is
divisible by $r = (\gamma_1\gamma_2)+(\gamma_1\gamma_2^{-1})-(\gamma_1)(\gamma_2)$, Lemma 3 implies
$\Psi(p')=0$.   Let $q'$ and $r'$ be the results of substituting
$w_i$ for $a_i$ in $p'$ and $r$ (respectively).   
Removing the root $cL$ from the 
forest for $\Psi(p)$ produces a forest for $\Psi(p')$ with lower
total potential.  Hence $\kappa(p') <\kappa(p)$ and, by induction,
$\Psi(q') =0$.  Furthermore, $r'$ has the form
$(\gamma_1'\gamma_2')+(\gamma_1'{\gamma_2'}^{-1})-(\gamma_1')(\gamma_2')$.  Since $r'$ divides $q-q'$
we have $\Psi(q)=0$.  

{\em Case 2: The resolution removes a crossing between two
components.}  In this case the components involved in the crossing
correspond to  loops $\gamma_1$ and $\gamma_2$, for which the resolution
produces $K_{\gamma_1\gamma_2}$ and $K_{\gamma_1\gamma_2^{-1}}$.  As in Case 1 we create
$p'$ by replacing $(\gamma_1)(\gamma_2)$ in $p$ with $(\gamma_1\gamma_2)+(\gamma_1\gamma_2^{-1})$.
The proof then proceeds by induction as before. 

{\em Case 3: The resolution removes a trivial circle.}  The trivial
circle corresponds to an appearance of $(1)$ in $p$.  Form $p'$ by
replacing that $(1)$ with the scalar $2$. Then create $q$ and $q'$ as
above, noting that $(1)-2$ divides both $p-p'$ and $q-q'$.  As above,
$\kappa(p') < \kappa(p)$, and it follows that $\Psi(q)=0$.
\end{proof}

We would now like to consider a more general sort of trace identity.
Let $\cS_n$ denote the group of permutations of the set
$\{a_1,\ldots,a_n\}$.  Let $\cS_m$ denote the group of permutations of
some subset $\{a_{i_1},\ldots,a_{i_m}\}$.  
Consider the group algebra $\bC\cS_m$.  By writing the elements of
$\cS_m$ in cycle notation, including trivial cycles, we obtain
expressions in $\bC[\cW]$. (Example 1, for instance.)  In fact, since
no inverses appear in these expressions, they can be regarded as
functions on the set of $m$-tuples of $2\times2$ matrices. If an
element of $\bC\cS_m$, regarded as such a function,  vanishes for
every assignment of $2\times2$ matrices we call it a {\em Procesi
identity} on $\cS_m$.   Note that a Procesi identity is clearly an
$\g$ trace identity, but that the converse is just as clearly false.

Using the group algebra to encode Procesi identities is useful for
the theorem we are about to prove, but there is a drawback as well.
Multiplication in $\bC\cS_m$ is not the same as multiplication in
$\bC[\cW]$.  If $p$ and $q$ are elements of $\bC\cS_m$ we denote their
product in the group algebra as $p\cdot q$, always assuming that $p$,
$q$ and $p \cdot q$ are written in cycle notation.   Note that $pq$
need not lie in $\bC\cS_m$, and that $p \cdot q$ may involve elements
of $\cW$ which do not appear in either $p$ or $q$.  Fortunately, the
skein module keeps track of how multiplication in
$\cS_m$ rearranges the elements of $\cW$. 

\begin{prop}
Let $p \in \bC\cS_m$.  If $\Psi(p)=0$ then $\Psi(\tau \cdot p) =0$ for every
$\tau \in \cS_n$.
\end{prop}

\begin{proof}
As an initial simplification assume that $\tau = (a_ia_j)$ with
$i<j$, and that $\cS_m$ permutes the set $\{a_1,\ldots,a_m\}$. There
are three cases, depending on the intersection of $\{a_1,\ldots,a_m\}$
and $\{a_i,a_j\}$.

{\em Case 1:} $m < i$.  As an element of $\bC[\cW]$, $\tau \cdot p$
factors into $(a_ia_j)p$. Hence $\Psi(p)=0$ implies $\Psi(\tau \cdot p)=0$.

{\em Case 2:} $i \leq m < j$.  Each term of $p$ contains a
cycle in which $a_i$ appears.  Assume that it is written
$(a_i\alpha)$ and write $\tau \cdot (a_i\alpha)$ as $(a_ja_i\alpha)$.
Fix a diagram for each term of $\Psi(p)$ with the property that it
traverses handles $1$ through $m$ exactly once and misses the others.
In a resolving forest for the standard resolution the skein relations
take place in neighborhoods away from the handles, and no trivial
circle runs once over a handle.  Therefore every diagram in the
forest meets the handles in precisely the same set of arcs, and we can
apply the operation shown in  Figure 5  
  \begin{figure}
    \mbox{}\hfill\epsfig{figure=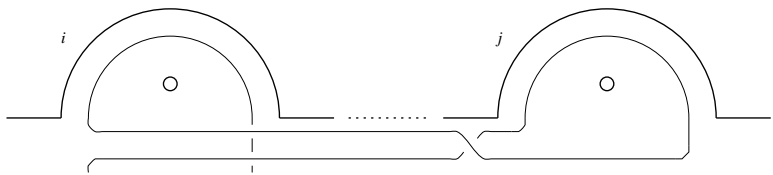,width=4in}\hfill\mbox{}
    \caption{Band sum of $\Psi(a_i\alpha)$ and $\Psi(a_j)$.}
  \end{figure}  
to the entire forest. Note that this changes the diagram for
$\Psi((a_i\alpha))$ into a diagram for 
$\Psi((a_ja_i\alpha))$, producing a resolution of
$\Psi(\tau \cdot p)$.  By Theorem 4, the resolution of $\Psi(p)$ is formally
zero in $\bC\cL_H$.  Since the resolution of $\Psi(\tau \cdot p)$ is
obtained by applying Figure 5 to each term, it must also be zero.

{\em Case 3:} $j \leq m$.
Each term of $p$ contains either a
cycle $(a_i\alpha a_j\beta)$ or a product of cycles
$(a_i\alpha)(a_j\beta)$.  The action of $\tau$  interchanges the
two possibilities.  Notice that the operation in Figure 6  
  \begin{figure}[b]
    \mbox{}\hfill\epsfig{figure=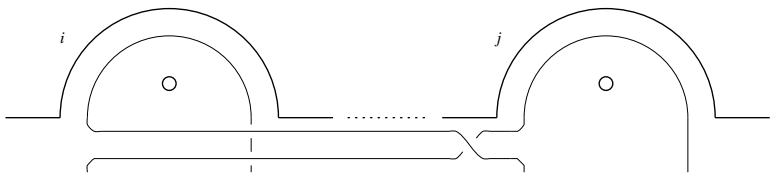,width=4in}\hfill\mbox{}
    \caption{Band relating $\Psi((a_i\alpha)(a_j\beta))$ and
             $\Psi(a_i\alpha a_j\beta)$.} 
  \end{figure}  
interchanges  the diagrams for $\Psi((a_i\alpha)(a_j\beta))$ and
$\Psi((a_i\alpha a_j\beta))$.  The proof then follows the resolving
argument of Case 2. 

Subject to our initial simplification, we now have $\Psi(\tau \cdot
p)=0$.  Retaining the assumption that $\tau = (a_ia_j)$, we next allow
$\cS_m$ to permute any set $\{a_{i_1},\ldots,a_{i_m}\}$. A
substitution converts this set into $\{a_1,\ldots,a_m\}$, but
preserves that fact that $\tau$ is a transposition.  Hence, by
Proposition 1, we again have $\Psi(\tau \cdot p)=0$.  Finally, since
any element of $\cS_n$ is a product of transpositions, $\Psi(\tau
\cdot p)=0$ for all $\tau \in \cS_n$.
\end{proof}

Our interest in Procesi identities stems from a classification 
theorem due independently to Procesi \cite{procesi} and Razmyslov
\cite{raz}.  Leron \cite{leron} is an excellent reference
for the proof. For the sake of completeness we include some
definitions taken from \cite[Chapter 4]{FH}.  A {\em
Young diagram} for $\cS_m$ is a collection of $m$ boxes arranged in left
justified rows of decreasing length.  A {\em Young tableau} is an
assignment of $a_{i_1},\ldots,a_{i_m}$ to the boxes.  Figure 7 is an
example of a Young diagram  for $S_{10}$ and a tableau using
$\{a_1,\ldots,a_{10}\}$.  
  \begin{figure}[t]
    \mbox{}\hfill\epsfig{figure=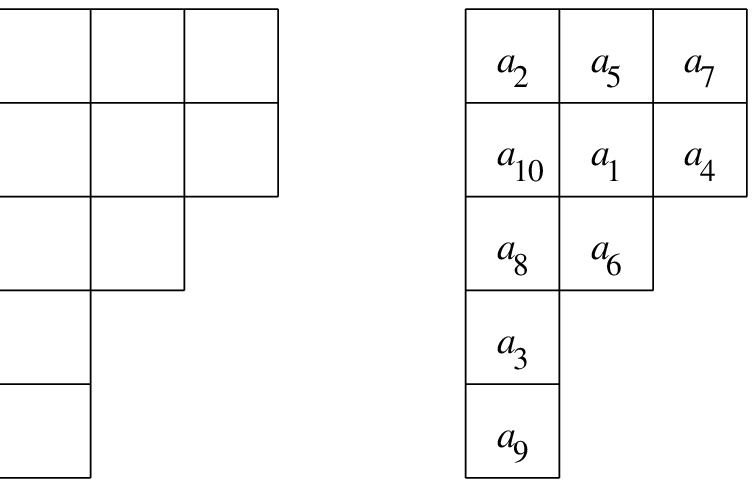,width=2.5in}\hfill\mbox{}
    \caption{}
  \end{figure}
Given a tableau $Y$ define $P_Y$ to be the subgroup of $\cS_m$ stabilizing 
the rows.  For the tableau in Figure 7 
\[P_Y \cong \cS_3 \times \cS_3\times
\cS_2.\]
Similarly, define $Q_Y$ to be the column stabilizer.
The Young symmetrizer corresponding to $Y$ is the element 
 \[ \left(\sum_{\sigma \in P_Y} \sigma\right) \cdot \left(\sum_{\tau \in
               Q_Y}\text{sgn}(\tau)\tau\right) \in \bC\cS_m. \]

\begin{theorem}{\rm (Procesi, Razmyslov)}
Procesi identities on a fixed $\cS_m$ constitute an ideal in
$\bC\cS_m$.  The ideal is generated by Young symmetrizers
corresponding to diagrams with at least three rows. 
\end{theorem}

\begin{lemma}
Let $Y$ be a Young tableau on $\{a_1,\ldots,a_m\}$, and assume 
that $a_m$ occupies the last box of a
row and column as shown in Figure 8.  Let $Y'$ be the tableau obtained
from $Y$ by removing the box containing $a_m$.  Using the notation
$r_{s+1}=c_{t+1}=a_m$, we can express $P_Y$ and $Q_Y$ as the following
disjoint unions:
\begin{enumerate}
\item ${\displaystyle P_Y=\bigcup_{i=1}^{s+1}(r_ia_m)\cdot P_{Y'}}$, and
\item ${\displaystyle Q_Y=\bigcup_{i=1}^{t+1}(c_ia_m)\cdot Q_{Y'}}$.
\end{enumerate}
\end{lemma}

  \begin{figure}[b]
    \mbox{}\hfill\epsfig{figure=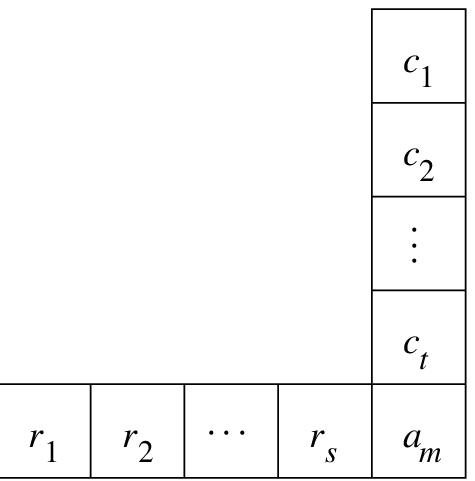,width=1.5in}\hfill\mbox{}
    \caption{}
  \end{figure}

\begin{proof}
Let $\lambda_1,\ldots,\lambda_x$ be the lengths of the rows of $Y'$.
A row stabilizer is a product of symmetric groups, so 
\[ |P_{Y'}|=\prod_{j=1}^x\lambda_j!,\quad\text{and} \]
\[ |P_{Y}|=(s+1)\prod_{j=1}^x\lambda_j!.\]
Each coset $(r_ia_m)\cdot P_{Y'}$ stabilizes the rows of $Y$ and, since each
contains the element $(r_ia_m)$, they are disjoint.  Counting elements
finishes the proof for $P_Y$.  The proof for $Q_Y$ is similar. 
\end{proof}

\begin{theorem}
If $p$ is a Procesi identity then $\Psi(p)=0$.
\end{theorem}

\begin{proof}
Implicit in the statement is the fact that $p$ is a Procesi identity
on some $\cS_m$. We proceed by induction on $m$.  By Theorem 5 and
Proposition 2 we may assume that $p$ is a Young
symmetrizer corresponding to a digram with at least three rows.  If
$m=3$ there is only one such diagram and $p$ is the result of
substituting $a_{i_1}$, $a_{i_2}$ and $a_{i_3}$ into Example 1.  By
Lemma 2 and Proposition 1, $\Psi(p)=0$.  

Now assume $m>3$.  Choose a diagram with at least three rows and a
tableau $Y$ satisfying the hypotheses of Lemma 3.   A symmetrizer
corresponding to any other tableau with the same diagram is obtained
from this one by by a substitution.  Therefore, by Proposition 1,
it suffices to consider only $Y$.  With notation as in Lemma 3, let
$p'$ be the symmetrizer corresponding to the tableau $Y'$.  We then have
\begin{equation*}\begin{split}
      p & =  \left(\sum_{\sigma \in P_Y} \sigma\right) \cdot
              \left(\sum_{\tau \in Q_Y}\text{sgn}(\tau)\tau\right)  \\
        & =  \left[\sum_{i=1}^{s+1}
               \left(\sum_{\sigma \in P_{Y'}} (r_ia_m)\cdot\sigma\right)
               \right]\cdot
               \left[\sum_{j=1}^{t+1}
                 \left(\sum_{\tau \in Q_{Y'}} \text{sgn}((c_ja_m)\cdot\tau)
                         (c_ja_m)\cdot\tau\right)
              \right] \\
        & =  \sum_{i,j} \text{sgn}(c_ja_m)(r_ia_m)\cdot(c_ja_m)\cdot
              \left(\sum_{\sigma \in P_{Y'}} \sigma\right)\cdot
              \left(\sum_{\tau\in Q_{Y'}}\text{sgn}(\tau)\tau\right)  \\
        & =  \sum_{i,j} \text{sgn}(c_ja_m)(r_ia_m)\cdot(c_ja_m)\cdot p'.
\end{split}\end{equation*}
By induction we have $\Psi(p)=0$.
\end{proof}

\section{The Coordinate Ring}

In Section 3 we developed conditions under which $\Psi$ carries an
$\g$ trace identity to zero.  In this section we will complete the
characterization of $\ker \Phi$, but to do so we must chose
coordinates on $X(H)$.
  
Not only did Vogt \cite{V}, Fricke \cite{F},  and Culler and
Shalen \cite{CS} apparently discover Theorem 1 independently, they all
arrived at the same set of generators.  Let $\gamma =
a_{i_1}\cdots a_{i_m}$ be an element of $G$ in which each $a_{i_j}$ is
distinct.  Following \cite{GM} we adopt the shorthand notation
$t_{i_1\cdots i_m}$ for the map $t_\gamma.$  The generating set in all
versions of Theorem 1 is $\cT=\{t_{i_1\cdots
i_m}\;|\;i_1<i_2<\cdots<i_m\}$.

Note that $\bC[\cT]$ becomes a subring of $\bC[\cW]$ by replacing
$t_{i_1\cdots i_m}$ with $(a_{i_1}\cdots a_{i_m})$, so $\Psi$ is well
defined on $\bC[\cT]$. Theorem 1 says that for every $p \in \bC[\cW]$
there exists $q \in \bC[\cT]$ such that $p$ and $q$ represent the same
element of $\cR(H)$. We will need a stronger result, for which we turn
to the combinatorial construction of $\cT$ in \cite{states}.

\begin{theorem} Let ${\mathcal K}$ be a set of knots containing exactly
one $K_\gamma$ for each $t_\gamma \in  \cT$.  Any link  $L \in \cL_H$ has
a resolving tree whose leaves are monomials in
$\bC[{\mathcal K}]$.
\end{theorem} 

\begin{cor}
For every $p \in \bC[\cW]$ there exists $q \in
\bC[\cT]$ such that $\Psi(p) = \Psi(q)$.  
\end{cor}

The main result of \cite{GM} is the construction of an ideal, $\cJ_H$,
which defines $X(H)$ in the coordinates $\cT_0 = \{t_{i_1\cdots
i_m} \in \cT \;|\; m \le 3\}$. The radical of
this ideal, $\sqrt{\cJ_H}$, is the ideal of trace identities in
$\bC[\cT_0]$.  The authors of \cite{GM} show that the 
trace identities in $\bC[\cT_0]$ generate those in $\bC[\cT]$, but we will
need a slightly stronger result.

\begin{lemma} {\rm (Compare \cite[Lemma 4.1.1]{GM}).} Choose distinct indices 
$i,j,k,m_1,\ldots,m_l$ and let $\alpha=m_1\cdots m_l$.  If
\begin{equation*}\begin{split} 
q  = & -2 t_{ijk\alpha} + t_{ik} t_j t_{\alpha} - t_i t_j t_{k\alpha}
        - t_j t_k t_{\alpha i} - t_{ik} t_{j\alpha} \\
     &  + t_{ij} t_{k\alpha} + t_{jk} t_{\alpha i} - t_{ikj}t_{\alpha}
        + t_i t_{jk\alpha} + t_j t_{k\alpha i} + t_k t_{\alpha ij}
\end{split}\end{equation*}
then $\Psi(q)=0$.  
\end{lemma}

\begin{proof}
Let $p_{xyz}$ denote the Procesi identity of Example 1 with the
substitutions $a_1=a_x$, $a_2=a_y$ and $a_3=a_z$.  Consider the polynomial 
\[p'=(a_1a_3)\cdot p_{234}-(a_3a_4)\cdot p_{124}-(a_1a_4)\cdot p_{123}\] 
as an element of $\bC[\cW]$.  By Theorem 5, $p'$ is a Procesi
identity, implying $\Psi(p')=0$.  Substituting
$a_1=a_i$, $a_2=a_j$, $a_3=a_k$ and $a_4=a_{m_1}\cdots a_{m_l}$ in
$p'$ we obtain $q$, so  Proposition 1 gives $\Psi(q)=0$.
\end{proof}

\begin{prop}
For every $q \in \bC[\cT]$ there exists $q_0 \in
\bC[\cT_0]$ such that $\Psi(q) = \Psi(q_0)$.  
\end{prop}

\begin{proof}
Let $q \in \bC[\cT]$.  Define $l$ to be the maximum length of a
subscript appearing in $q$ and let $m$ be the number of maximum
length subscripts.  We say that the complexity of $q$ is the ordered
pair $(l,m)$.  The proof is by induction on complexity, ordered
lexicographically.  If $l \le 3$ then $q$ can be converted to $q_0 \in
\bC[\cT_0]$ by repeated application of the identity in Example
1 (perhaps with a substitution of indices).  The difference between
any pair of successive stages is divisible by a Procesi identity, so
Theorem 6 implies $\Psi(q)=\Psi(q_0)$.

If  $l>3$ then $q$ contains some $t_{ijk\alpha}$ in which $ijk\alpha$
is a maximum length subscript.  The identity of
Lemma 4 allows us to replace $t_{ijk\alpha}$ with an expression
involving only shorter subscripts.  The result is a new polynomial
$q'$ with lower complexity and $\Psi(q)=\Psi(q')$.  
\end{proof}

We now state the main result of \cite{GM}.  Define 
\begin{align*}
M_{ii} & = t_{i}^2-4, \qquad \text{and} \\
M_{ij} = M_{ji} &= 2t_{ij}-t_it_j, \qquad \text{if} \quad i<j.
\end{align*}

\begin{theorem} {\rm (Gonz\'{a}lez-Acu\~{n}a--Montesinos)} $X(M)$ is
the zero set of the ideal $\cJ_H$ in $\bC[\cT_0]$ generated by the
following polynomials.
\begin{align*}
\begin{split}
q_1 & =  t_i^2 + t_j^2 + t_k^2 + t_{ij}^2 +t_{ik}^2 + t_{jk}^2
       +t_{ijk}^2  + t_{ij} t_{ik} t_{jk} + t_{ijk}t_it_jt_k\\
    &\quad -t_{ijk}t_i t_{jk}-t_{ijk}t_j t_{ik} -t_{ijk}t_k t_{ij}
       - t_i t_j t_{ij}- t_i t_k t_{ik} - t_j t_k t_{jk}-4,\\
    &\quad\mbox{in which $i$, $j$ and $k$ are distinct.}
\end{split} \\
q_2 &= \begin{vmatrix}
      M_{11} & M_{12} & M_{1i} & M_{1j} \\
      M_{21} & M_{22} & M_{2i} & M_{2j} \\
      M_{i1} & M_{i2} & M_{ii} & M_{ij} \\ 
      M_{j1} & M_{j2} & M_{ji} & M_{jj} 
      \end{vmatrix}, \quad\mbox{with $2 < i < j \leq n$.} \\
q_3 &= \begin{vmatrix}
      M_{11} & M_{12} & M_{13} & M_{1i} \\
      M_{21} & M_{22} & M_{23} & M_{2i} \\
      M_{31} & M_{32} & M_{33} & M_{3i} \\ 
      M_{j1} & M_{j2} & M_{j3} & M_{ji} 
      \end{vmatrix}, \quad\mbox{with $3 < i < j \leq n$.} \\
q_4 &= (t_{123}-t_{132})
        (2t_{ijk}+t_i t_j t_k-t_i t_{jk}-t_j t_{ik}-t_k t_{ij}) 
    - \begin{vmatrix}
      t_1 & t_{1i} & t_{1j} & t_{1k} \\
      t_2 & t_{2i} & t_{2j} & t_{2k} \\
      t_3 & t_{3i} & t_{3j} & t_{3k} \\
       2  & t_i    &    t_j & t_k       
      \end{vmatrix},\\
    &\quad\mbox{in which $1\leq i<j<k\leq n$ and $t_{mm}$ denotes $t_m^2-2$.}
\end{align*}
\end{theorem}

The next step is to show that all of these polynomials lie in $\ker
\Psi$.  The proofs involved in this are closely modeled on those in
\cite{GM}.  Our contribution is the observation that $q_2$, $q_3$
and $q_4$ lift to Procesi identities, and that $q_1$ follows directly
from a resolving tree.  We will introduce further notation from
\cite{GM} as it becomes necessary.

\begin{lemma} $\Psi(q_1)=0$.
\end{lemma}

\begin{proof}
Example 2, Lemma 1 and Proposition 1.
\end{proof}

Let $A_1,\ldots,A_4$ be $2\times2$ matrices with
\[A_i = \begin{pmatrix}
        \alpha_i & \beta_i \\
        \gamma_i & \delta_i
        \end{pmatrix}.\]
Define 
\[M(A_i) = \begin{pmatrix}
           \alpha_1 & \beta_1 & \gamma_1 & \delta_1 \\
           \alpha_2 & \beta_2 & \gamma_2 & \delta_2 \\
           \alpha_3 & \beta_3 & \gamma_3 & \delta_3 \\
           \alpha_4 & \beta_4 & \gamma_4 & \delta_4
           \end{pmatrix},\;
  J = \begin{pmatrix}
      1 & 0 & 0 & 0 \\
      0 & 0 & 1 & 0 \\
      0 & 1 & 0 & 0 \\
      0 & 0 & 0 & 1 
      \end{pmatrix}, \;\text{and} \;
  J^*= \begin{pmatrix}
      0 & 0 & 0 & 1 \\
      0 & -1 & 0 & 0 \\
      0 & 0 & -1 & 0 \\
      1 & 0 & 0 & 0 
      \end{pmatrix}.\]
In general, $(c_{ij})$ will denote a $4\times4$ matrix and
$|c_{ij}|$ its determinant.  We use $A^t$ to denote the transpose
of $A$.  

\begin{lemma} {\rm (Compare \cite[Lemma 4.6]{GM}).}  The polynomial
$p = |2(x_iy_j) - (x_i)(y_j)|$ is a Procesi identity on the
symmetric group permuting $\{x_1,x_2,x_3,x_4,y_1,y_2,y_3,y_4\}$.
\end{lemma}

\begin{proof}
Clearly $p$ is an element of the group algebra over the permutations
of $\{x_1,x_2,x_3,x_4,y_1,y_2,y_3,y_4\}$.  To see that $p$ is an
identity, assign matrices $A_i$ and
$B_j$ to each $x_i$ and $y_j$.  Direct calculation shows that
\[ (\tr(A_iB_j)) = M(A_i)JM(B_j)^t, \]
and 
\[ (\tr(A_iB_j)-\tr(A_i)\tr(B_j)) = -M(A_i)JJ^*M(B_j)^t.\]
Since $|I-J^*|=0$ we have
\[|2\;\tr(A_iB_j) - \tr(A_i)\tr(B_j)| =  |M(A_i)J(I-J^*)M(B_j)^t|= 0.\]
\end{proof}

\begin{lemma} {\rm (Compare \cite[Corollary 4.12]{GM}).} The polynomial
\begin{equation*}\begin{split}
p &=[(x_1x_2x_3)-(x_1x_3x_2)] \\
    & \quad \times 
 [2(y_1y_2y_3)+(y_1)(y_2)(y_3)-(y_1)(y_2y_3)-(y_2)(y_1y_3)-(y_3)(y_1y_2)]\\
    & \quad -\begin{vmatrix}
      (x_1) & (x_1y_1) & (x_1y_2) & (x_1y_3) \\
      (x_2) & (x_2y_1) & (x_2y_2) & (x_2y_3) \\      
      (x_3) & (x_3y_1) & (x_3y_2) & (x_3y_3) \\
        2   &  (y_1)   &  (y_2)   &  (y_3)   
      \end{vmatrix} 
\end{split}\end{equation*}
is a Procesi identity on $\{x_1,x_2,x_3,y_1,y_2,y_3\}$. 
\end{lemma}

\begin{proof}(Compare \cite[Lemma 4.10 and Proposition 4.11]{GM}).
For any $2\times2$ matrices $A_1,\ldots,A_4,B_1,\ldots,B_4$ we have
\begin{equation*}\begin{split}
|\tr(A_iB_j)| & =  |M(A_i)JM(B_j)| \\    
              & =  |J||M(A_i)||M(B_j)| \\
              & =  -|M(A_i)||M(B_j)|.
\end{split}\end{equation*}
If $A_4=I$ then, by direct calculation, 
\[ |M(A_i)| = \tr(A_1A_2A_3)-\tr(A_1A_3A_2). \]
If $B_4=I$ as well then
\[ |\tr(A_iB_j)|= -[\tr(A_1A_2A_3)-\tr(A_1A_3A_2)]
                   [\tr(B_1B_2B_3)-\tr(B_1B_3B_2)].  \]
Changing columns in $|\tr(A_iB_j)|$ and applying the identity of
Example 1, we see that $p$ vanishes for an arbitrary assignment of $2\times 2$ matrices.  As in
Lemma 6, it is clearly an element of the appropriate group
algebra.
\end{proof}
 
\begin{lemma} {\rm (Compare \cite[Proposition 4.8]{GM}).} $\Psi(q_2) = 0$.
\end{lemma}

\begin{proof}
Create $p \in \bC[\cW]$ by specializing the Procesi identity of Lemma
6 at $x_1=y_1 = a_1$, $x_2=y_2 = a_2$, $x_3=y_3 = a_i$ and $x_4=y_4 =
a_j$.  Theorem 6 and Proposition 1 imply $\Psi(p)=0$.  Note that $p$
and $q_2$ are determinants of matrices that differ only along their
diagonals.  The differences between diagonal terms are of the form
\[ 2(a_m^2)-(a_m)^2-t_m^2+4,\]
which can be rewritten as
\[ 2(a_m^2) -2(a_m)^2+2(1)-2(1)+4.\]
Hence, $q_2$ may be obtained from $p$ by a finite sequence of
substitutions of the form $(a_m^2) = (a_m)^2-(1)$ or $(1)=2$.  Each
step involves a pair of polynomials whose difference is divisible
either by $(a_ma_m)+(a_ma_m^{-1})-(a_m)(a_m^{-1})$ or by $(1)-2$, so
Lemma 2 implies $\Psi(q_2)=\Psi(p)=0$.
\end{proof}

\begin{lemma} {\rm (Compare \cite[Proposition 4.9]{GM}).} $\Psi(q_3) = 0$.
\end{lemma}

\begin{proof}
Specialize the identity of Lemma 6 at 
$x_1=y_1 = a_1$, $x_2=y_2 = a_2$, $x_3=y_3 = a_3$, $x_4=a_i$ and $y_4 = a_j$.
Then proceed as in Lemma 8.
\end{proof}

\begin{lemma} {\rm (Compare \cite[Corollary 4.12]{GM}).} $\Psi(q_4)=0$. 
\end{lemma}

\begin{proof}Specialize the Procesi
identity of Lemma 7 at $y_1=a_i$, $y_2=a_j$, $y_3=a_k$ and $x_m = a_m$
for $m =1,2,3$.  If $i>3$ this is precisely $q$.  If not, then
proceed as in Lemmas 8 and 9, using the fact that $t_{mm}$ denotes
$t_m^2-2$.
\end{proof}

This is enough to prove our claims about $\Phi : V(M) \rightarrow
\cR(M)$, but we may as well consider an arbitrary, compact, orientable
3-manifold.  Let $M$ be the result of adding 2-handles to $H$ along
curves $\{c_1,\ldots,c_m\}$ in $\partial H$.  Choose words $w_i$ in
$\pi_1(H)$ so that, as a loop, each $w_i$ is freely homotopic to some
orientation of $c_i$.  For each $i$ and $j$ form the polynomial
$p_{ij} = (w_ia_j)-(a_j) \in \cR(H)$.  Using the obvious
identification $\cR(H) \cong \bC[\cT_0]/\sqrt{\cJ_H}$, create an ideal
$\cJ_M$ in $\bC[\cT_0]$ generated by $\cJ_H \cup \{p_{ij}\}$.

\begin{theorem}  
{\rm (Gonz\'{a}lez-Acu\~{n}a--Montesinos)} $X(M)$ is the zero set of
$\cJ_M$ in $\bC[\cT_0]$.
\end{theorem}

It follows immediately that $\cR(M) = \bC[\cT_0]/\sqrt{\cJ_M}$.  We
know that $\Phi$ maps $V(M)$ onto $\cR(M)$ and it is clear that $\Psi$
maps $\bC[\cW]$ onto $V(H)$, and hence onto $V(M)$.  Using these maps,
we can now see how $V(M)$ compares to $\cR(M)$.  From now on, consider
$\Psi$ to be the restriction to $\bC[\cT_0]$. 

\begin{prop} 
$\cJ_M \subset \ker \Psi \subset \sqrt{\cJ_M}$
\end{prop}

\begin{proof}
That $\cJ_H$ lies in $\ker \Psi$ is the content of Lemmas 5, 8, 9, and
10.  To see that $\Psi(p_{ij})=0$, construct a knot $K_{a_j}$ for each
generator $a_j$.  For each $i$ and $j$, there is a band sum $c_i \#_b
K_{a_j}$ producing a knot $K_{w_ia_j}$.  Since $c_i\#_b K_{a_j} \cong
K_{a_j}$  in $M$, we have $\Psi(p_{ij})=K_{a_j}-K_{w_ia_j}=0$. 

For the second containment, note that Corollary 1 and Proposition 3
imply that $\Psi|_{\bC[\cT_0]}$ is still onto.  It should now be clear
that
\[\bC[\cT_0] \stackrel{\Psi}{\rightarrow} V(M)
              \stackrel{\Phi}{\rightarrow} \cR(M) 
              \cong \bC[\cT_0]/\sqrt{\cJ_M}\]
is the canonical projection.  
\end{proof}

There are various equivalent ways of phrasing the immediate
implication of Proposition 4.

\begin{theorem} Let $M$ be a compact orientable 3-manifold with
$\Phi$, $\Psi$, $\cJ_M$, and $\cT_0$ defined  as above.  Denote the
ideal of nilpotents in $V(M)$ by $\sqrt{0}$. 
\begin{enumerate}
\item $X(M)$ is the zero set of $\ker \Psi$ in $\bC[\cT_0]$. 
\item $\sqrt{\ker \Psi} = \sqrt{\cJ_M}$.
\item $\ker \Phi = \sqrt{0}$.
\item $\Phi$ induces an isomorphism $\widehat{\Phi}:V(M)/\sqrt{0}
\rightarrow \cR(M)$.
\item $\Psi$ induces an isomorphism $\widehat{\Psi}:
\bC[\cT_0]/\sqrt{\cJ_M} \rightarrow V(M)/\sqrt{0}$.
\item Under the identification of $\cR(M)$ with
$\bC[\cT_0]/\sqrt{\cJ_M}$, the maps $\widehat{\Psi}$ and
$\widehat{\Phi}$ are inverses. 
\end{enumerate}
\end{theorem}

\begin{proof}\mbox{}
\begin{enumerate}
\item This is immediate from Proposition 4 and the fact the $X(M)$ is
the zero set of both $\cJ_M$ and $\sqrt{\cJ_M}$.
\item Nullstellensatz.
\item Since $\cR$ cannot, by definition, contain a non-zero nilpotent
element, $\Phi(\sqrt{0}) = 0$.  Suppose now that $\Phi(\alpha)=0$, and
write $\alpha$ as $\Psi(\beta)$.  We have seen that
\[\bC[\cT_0] \stackrel{\Psi}{\rightarrow} V(M)
\stackrel{\Phi}{\rightarrow} \cR(M) \cong \bC[\cT_0]/\sqrt{\cJ_M}\] is
the canonical projection.  Hence, $\beta \in \sqrt{\cJ_M}$.  It
follows from Theorem 10 part (2) that $\Psi(\beta^n) = 0$ for some $n$, meaning
$\alpha$ is nilpotent.
\item Theorem 3 and Theorem 10 part (3).
\item The composition 
\[\bC[\cT_0] \stackrel{\Psi}{\rightarrow} V(M)
\stackrel{\pi}{\rightarrow} V(M)/\sqrt{0}
\stackrel{\widehat{\Phi}}{\rightarrow} \cR(M) \cong
\bC[\cT_0]/\sqrt{\cJ_M}.\] is also the canonical projection (here
$\pi$ is projection as well).  Hence, $\ker \pi \circ \Phi =
\sqrt{\cJ_M}$.
\item It is easy to see that both
\[\bC[\cT_0]/\sqrt{\cJ_M}
 \stackrel{\widehat{\Psi}}{\rightarrow} V(M)/\sqrt{0}
\stackrel{\widehat{\Phi}}{\rightarrow} \cR(M) \cong
\bC[\cT_0]/\sqrt{\cJ_M}\]
and
\[ V(M)/\sqrt{0}\stackrel{\widehat{\Phi}}{\rightarrow} \cR(M) \cong
\bC[\cT_0]/\sqrt{\cJ_M}\stackrel{\widehat{\Psi}}{\rightarrow}
V(M)/\sqrt{0}\]
are the identity.  
\end{enumerate}
\end{proof}

We conclude with an application.  The author would like to thank
Charles Frohman for suggesting that this result might follow quickly,
Victor Camillo for encouraging us to disregard nilpotents, and
Bernadette Mullins for pointing out the result from ring theory used in
the proof.  Recall that a 3-manifold is {\em small} if it contains no
incompressible, non-boundary parallel surface.

\begin{theorem} {\rm (Compare \cite[Corollary 1]{reps})}. 
If $M$ is small then $\dim V(M) < \infty$. 
\end{theorem} 

\begin{proof}
Suppose that $X(M)$ has positive dimension.  If follows that some
component of $X(M)$ contains a curve whose smooth projective
resolution contains an ideal point.  From \cite[2.2.1]{CS} we then
have a non-trivial splitting of $\pi_1(M)$, meaning $M$ is not small.
Hence, $X(M)$ consists of a finite set of points and $\cR(M)$ is
finite dimensional as a vector space.  It is a standard result of
commutative algebra that an ideal in a Noetherian ring contains some
power of its radical.  Thus, from Theorem 10 part (2), we obtain
\[\left(\sqrt{\cJ_M}\right)^n \subset \ker \Psi \subset \sqrt{\cJ_M}\]
for some $n$.  Since $\cR(M) \cong \bC[\cT_0]/\sqrt{\cJ_M}$, it is a
simple exercise to show that $\bC[\cT_0]/\left(\sqrt{\cJ_M}\right)^n$
is finite dimensional.  The result now follows from the fact that
$V(M) \cong \bC[\cT_0]/\ker \Psi$, which in turn is the homomorphic
image of $\bC[\cT_0]/\left(\sqrt{\cJ_M}\right)^n$.
\end{proof}

\end{document}